\documentclass[17pt,envcountsame]{extarticle}

\usepackage{amsmath,amssymb,amsthm}
\usepackage[utf8]{inputenc}
\usepackage{authblk}

\usepackage {graphicx}
\newtheorem{theorem}{Theorem}[section]
\newtheorem{lemma}[theorem]{Lemma}
\newtheorem{remark}[theorem]{Remark}
\newtheorem {proposition}[theorem]{Proposition}

\newtheorem {corollary}[theorem]{Corollary}
\newtheorem{problem}[theorem]{Problem}
\def \RR{\mathbb R}
\def \NN{\mathbb N}
\def\({\left(}  \def\){\right)}
\def\[{\left[}  \def\]{\right]}

\def \beq {\begin {equation}}
\def \eeq {\end{equation}}
\def \OL {\overline}

\def\NN {\rm {\bf N}}

\begin {document}

\title {Polarization problem on a higher-dimensional sphere for a simplex}

\author {Sergiy Borodachov}

\affil {Towson University, Towson, MD, USA 21252}

\maketitle

\abstract {We study the problem of maximizing the minimal value over the sphere $S^{d-1}\subset \RR^d$ of the potential generated by a configuration of $d+1$ points on $S^{d-1}$ (the maximal discrete polarization problem). The points interact via the potential given by a function $f$ of the Euclidean distance squared, where $f:[0,4]\to (-\infty,\infty]$ is continuous (in the extended sense) and decreasing on $[0,4]$ and finite and convex on $(0,4]$ with a concave or convex derivative $f'$. We prove that the configuration of the vertices of a regular $d$-simplex inscribed in $S^{d-1}$ is optimal. This result is new for $d>3$ (certain special cases for $d=2$ and $d=3$ are also new). As a byproduct, we find a simpler proof for the known optimal covering property of the vertices of a regular $d$-simplex inscribed in $S^{d-1}$.}

\bigskip

\noindent {Keywords: generalized Chebyshev constant, maximal polarization, potential, sphere, simplex, optimal covering problem}

\section {Introduction}

The problem of maximal discrete polarization on a given compact set in the Euclidean space $\RR^d$ was studied by many authors in recent years, see, in particular, papers \cite {Sto1975circle,Amb2009,AmbBalErd2013,ErdSaf2013,HarKenSaf2013,BorBos2014,Su2014,BorHarRez2018,RezSafVol2018,HarPetSafsub}. 
It requires placing a given number $N$ of equal point charges of the same sign on a given conductor so that the minimal value of the total potential on the conductor is as large as possible. Of a particular interest are solutions to this problem when the conductor is a sphere, where very few exact solutions are currently known.

Let $S^{d-1}:=\{(x_1,\ldots,x_d)\in \RR^d : x_1^2+\ldots+x_d^2=1\}$, $d\geq 2$, denote the unit sphere in the Euclidean space $\RR^d$. For a lower semi-continuous function $f:[0,4]\to (-\infty,\infty]$ (called the potential function) and a $(d+1)$-point configuration $\omega_d=\{{\bf v}_0,{\bf v}_1,\ldots,{\bf v}_d\}$ on $S^{d-1}$, denote
\begin {equation}\label {potential}
p_f({\bf x},\omega_d):=\sum\limits_{i=0}^{d}f\(\left|{\bf x}-{\bf v}_i\right|^2\),\ \ \ {\bf x}\in S^{d-1},
\end {equation}
where $\left|\ \!\cdot\ \!\right|$ denotes the Euclidean norm in $\RR^d$, and let
\begin {equation}\label {k}
P_{f}(\omega_d,S^{d-1}):=\min\limits_{{\bf x}\in S^{d-1}}p_f({\bf x},\omega_d).
\end {equation}
We consider the following problem.
\begin {problem}\label {P1}
{\rm
Find the quantity
\begin {equation}\label {p2}
\mathcal P_f(S^{d-1}):=\sup\{P_f(\omega_d,S^{d-1}) : \omega_d\subset S^{d-1},\  \# \omega_d=d+1\}
\end {equation}
and $(d+1)$-point configurations on $S^{d-1}$ that attain the supremum on the right-hand side of \eqref {p2}, which we will call optimal configurations.
}
\end {problem}
The constant $\mathcal P_f(S^{d-1})$ is known as the $(d+1)$-point $f$-polarization of $S^{d-1}$ or the $(d+1)$-th Chebyshev $f$-constant of $S^{d-1}$. In the general polarization problem instead of the sphere one considers arbitrary infinite compact set $A\subset \RR^d$ and point configurations $\omega_N\subset A$ of arbitrary fixed cardinality $N\geq 1$.

An important special case of the potential function $f$ in \eqref {potential} is
$$
f_s(t):=\begin {cases}
t^{-s/2}, & s>0,\cr
\frac {1}{2}\ln \frac {1}{t}, & s=\log,\cr
-t^{-s/2}, & s<0,\cr
\end {cases}
$$
which corresponds to the case of the Riesz $s$-kernel $f_s(\left|{\bf x}-{\bf y}\right|^2)=\left|{\bf x}-{\bf y}\right|^{-s}$ for $s>0$ and $f_s(\left|{\bf x}-{\bf y}\right|^2)=-\left|{\bf x}-{\bf y}\right|^{-s}$ for $s<0$ as well as the logarithmic kernel $f_{\log}(\left|{\bf x}-{\bf y}\right|^2)=\ln \frac {1}{\left|{\bf x}-{\bf y}\right|}$.
The exact solution to the maximal polarization problem on the unit circle $S^1$ is known for every $N\geq 1$ and the potential functions $f_s$ for the following ranges of $s$. The case $N=3$ was settled by Stolarsky  \cite {Sto1975circle} for $-2<s<0$, by Nikolov and Rafailov \cite {NikRaf2011} for $s>0$ and $s<-2$ (the case $s=-2$ is trivial), and for $s=\log$ by Hardin, Kendall, and Saff \cite{HarKenSaf2013}. For the case of arbitrary $N\geq 4$, see works by Ambrus \cite {Amb2009} and Ambrus, Ball, and Erd\'elyi \cite {AmbBalErd2013} ($s=2$), Erd\'elyi and Saff \cite {ErdSaf2013} ($s=4$), and Hardin, Kendall, and Saff \cite{HarKenSaf2013} (arbitrary $s\geq -1$, $s\neq 0$, and $s=\log$). The configuration of the vertices of a regular $N$-gon inscribed in $S^1$ was shown to be optimal in these papers. In fact, in \cite {HarKenSaf2013}, the optimality of the vertices of a regular $N$-gon inscribed in $S^1$ was established for any kernel of the form $f(\ell ({\bf x},{\bf y}))$, where $\ell({\bf x},{\bf y})$ is the geodesic distance between points ${\bf x}$ and ${\bf y}$ on $S^1$ and $f$ is non-increasing and convex on $(0,\pi]$ and continuous (in the extended sense) at $0$.

On a higher-dimensional sphere, some cases are known to have a simple proof (see, e.g.,  the book \cite [Section 14.2]{BorHarSafbook}). Namely, for $2\leq N\leq d$ points on the sphere $S^{d-1}$, the solution to the maximal polarization problem is known for any non-increasing and convex potential function $f$. It is any $N$-point configuration on $S^{d-1}$ with its center of mass at the origin. Furthermore, for $s=-2$, and arbitrary cardinality $N\geq 2$, the solution is any $N$-point configuration with center of mass at the origin.

One non-trivial result is known for the sphere $S^2\subset \RR^3$, where Problem \ref {P1} was solved for $N=4$ points and potential functions of the form
\begin {equation}\label {Rafailov}
f(t)={\rm sgn}(s)(t+C)^{-s/2}, \ \ \ C\geq 0, 
\end {equation}
with $s>-2$, $s\neq 0$, see the thesis by Su \cite {Su2014} and the result of Nikolov and Rafailov \cite {NikRaf2013} finding quantity \eqref {k} for a regular simplex and potential functions of form \eqref {Rafailov}.
The optimality of the configuration of the vertices of a regular simplex inscribed in $S^2$ was proved. 

Throughout the rest of the paper, $\omega_d^\ast:=\{{\bf x}_0,{\bf x}_1,\ldots,{\bf x}_d\}$ will denote the set of vertices of a regular $d$-simplex inscribed in $S^{d-1}$. The result actually proved in \cite {Su2014} asserts that for a potential function $f$ non-increasing and convex on $(0,4]$, if the absolute minimum over $S^{2}$ of the potential $p_f({\bf x},\omega_3^\ast)$ is achieved at points of the set $-\omega_3^\ast$  (antipodes of the points from $\omega_3^\ast$), then $\omega_3^\ast$ is optimal for the maximal polarization problem for $N=4$ points on $S^2$. However, up to this point, the absolute minimum of $p_f({\bf x},\omega_d^\ast)$ was shown to be achieved at points of $-\omega_d^\ast$ only for $f$ of form \eqref {Rafailov}, see \cite {NikRaf2013} and references therein. Furthermore, one can construct potential functions $f$ non-increasing and convex on $(0,4]$ such that the potential $p_f({\bf x},\omega_d^\ast)$, $d\geq 2$, does not achieve its absolute minimum over $S^{d-1}$ at points of $-\omega_d^\ast$ (see Corollary \ref {minimum'} below).

In the current paper, we obtain the solution to Problem \ref {P1} for $N=d+1$ points on $S^{d-1}$ (in any dimension $d$) and convex and non-increasing potential functions $f$ with a concave or convex derivative $f'$ (see Theorems \ref {optimalsimplex} and \ref {optimalsimplex'}). We also show that the absolute minimum over $S^{d-1}$ of the potential $p_f({\bf x},\omega_d^\ast)$ is achieved at points of the set $-\omega_d^\ast$ when $f$ has a concave derivative on $(0,4]$ and at points of $\omega_d^\ast$ when $f$ is finite at $t=0$ and has a convex derivative on $(0,4]$, see Theorem \ref {minimum} and Corollary \ref {minimum'} (a more detailed review on this problem is given in Section \ref {mainresults}).

The assumptions in Theorems \ref {optimalsimplex} and \ref {optimalsimplex'} hold for potential functions with $f'\leq 0$, $f''\geq 0$, and $f'''\leq 0$ or $f'''\geq 0$ on the interval $(0,4]$ (continuous in the extended sense at $0$). In particular, Theorem \ref {optimalsimplex} is valid for any completely monotone potential function on $(0,4]$ defined at $0$ by its limit value. Recall that an infinitely differentiable function $f$ is called completely monotone on an interval $I$ if $(-1)^k f^{(k)}\geq 0$ on $I$ for every $k\in \NN\cup \{0\}$. Such is the potential function $f_s$ defining the Riesz $s$-kernel for $s>0$ and, after adding an appropriate positive constant, the Riesz $s$-kernel for $-2<s<0$ and the logarithmic kernel.
Another important example of a kernel defined by a completely monotone potential function  is the Gaussian kernel $f(\left|{\bf x}-{\bf y}\right|^2)={\rm exp}(-\sigma \left|{\bf x}-{\bf y}\right|^2)$, where $\sigma$ is a positive constant. 

We remark that Gaussian kernels cannot be represented as a non-increasing and convex function of the geodesic distance $\ell({\bf x},{\bf y})$. At the same time certain kernels of the form $f(\ell ({\bf x},{\bf y}))$, where $f$ is  decreasing and convex, cannot be given as a convex function of $\left|{\bf x}-{\bf y}\right|^2$ (take, say $f(t)=-t$).

For sets other than the sphere, certain exact results on the maximal polarization are also known. For a $d$-dimensional ball and $-2<s\leq d-2$, $s\neq 0$, or $s=\log$, the $N$-point configuration with all its points located at the center of the ball is optimal for every $N\geq 1$, see \cite {ErdSaf2013}.  
In the case of the logarithmic potential ($s=\log$), the maximal polarization problem on an infinite compact set in the complex plane is solved for every $N$ by the zeros of the restricted Chebyshev polynomial for the set $A$ (see, e.g., \cite [Section 14.2]{BorHarSafbook}).

A number of recent works also deals with asymptotic behavior of the maximal polarization problem on various classes of compact rectifiable sets in $\RR^d$ (see the works by Erd\'elyi and Saff \cite {ErdSaf2013}, the author and Bosuwan \cite {BorBos2014}, the author, Hardin, Reznikov, and Saff \cite {BorHarRez2018}, and Reznikov, Saff, and Volberg \cite {RezSafVol2018}) as well as of the unconstrained maximal polarization problem, where configurations are allowed to lie anywhere in $\RR^d$ while the minimum of their potential is taken over the set $A$ (see the work by Hardin, Petrache, and Saff~\cite {HarPetSafsub}). The continuos version of the maximal polarization problem; i.e., when the potential of a finite configuration in \eqref {potential} is replaced by the potential of a Borel probability measure and one searches for a measure with the largest minimum of the potential on the conductor, was, in particular, considered Ohtsuka \cite {Oht1967}, Farkas and R\'ev\'esz \cite{FarRev2006}, Farkas and Nagy \cite {FarNag2008}, and Simanek \cite {Sim2016}. More information and references on the maximal polarization problem can be found, for example, in \cite [Chapter 14]{BorHarSafbook}.

\section {Main results}\label {mainresults}

The solution to Problem \ref {P1} is given by the two theorems below. 
\begin {theorem}\label {optimalsimplex}
Let $d\geq 2$ and $f:[0,4]\to(-\infty,\infty]$ be finite, non-increasing, and convex on $(0,4]$, differentiable in $(0,4)$ with a concave derivative $f'$ on $(0,4)$ such that $\lim\limits_{t\to 0^+}f(t)=f(0)$. Then
\begin {equation}\label {optimal}
P_f(\omega_d,S^{d-1})\leq P_f(\omega_d^\ast,S^{d-1})
\end {equation}
for every $(d+1)$-point configuration $\omega_d\subset S^{d-1}$. Furthermore,
$$
\mathcal P_f(S^{d-1})=P_f(\omega_d^\ast,S^{d-1})=f(4)+d\!\cdot\! f\(2-\frac {2}{d}\).
$$

If the convexity of $f$ is strict on $(0,4]$, then equality in \eqref {optimal} holds if and only if $\omega_d$ is the set of vertices of a regular $d$-simplex inscribed in $S^{d-1}$.
\end {theorem} 
We remark that part of the case $d=2$ of Theorem \ref {optimalsimplex} follows from the results in \cite {Sto1975circle,Amb2009,NikRaf2011,AmbBalErd2013,ErdSaf2013,HarKenSaf2013} and the case $d=3$ for potential functions of form \eqref {Rafailov}
follows by combining the result of \cite {Su2014} (which is the case $d=3$ of Lemmas \ref {semi} and \ref {polarizat}) with one of the results from \cite {NikRaf2013}. The case $d=3$ for general potentials follows by combining the result from \cite {Su2014} with the assertion of Theorem \ref {minimum} below. 


We also show the optimality of $\omega_d^\ast$ in the following case.
\begin {theorem}\label {optimalsimplex'}
Let $d\geq 2$ and $f:[0,4]\to(-\infty,\infty)$ be non-increasing and convex on $(0,4]$, differentiable in $(0,4)$ with a convex derivative $f'$ on $(0,4)$ such that $\lim\limits_{t\to 0^+}f(t)=f(0)$. Then
\begin {equation}\label {optimal'}
P_f(\omega_d,S^{d-1})\leq P_f(\omega_d^\ast,S^{d-1})
\end {equation}
for every $(d+1)$-point configuration $\omega_d\subset S^{d-1}$. Furthermore,
$$
\mathcal P_f(S^{d-1})=P_f(\omega_d^\ast,S^{d-1})=f(0)+d\!\cdot\! f\(2+\frac {2}{d}\).
$$

If the convexity of $f$ is strict on $( 0,4]$, then equality in \eqref {optimal'} holds if and only if $\omega_d$ is the set of vertices of a regular $d$-simplex inscribed in $S^{d-1}$.
\end {theorem} 
\begin {remark}
{\rm
The assumptions about the potential function $f$ in Theorem \ref {optimalsimplex'} are equivalent to the assumption that $f$ is continuously differentiable, non-increasing, and convex on $[0,4]$ with a convex derivative $f'$ on $[0,4]$.
}
\end {remark}

An important ingredient of the proof of Theorems \ref {optimalsimplex} and \ref {optimalsimplex'} is finding the absolute minimum over $S^{d-1}$ of the potential $p_f({\bf x},\omega_d^\ast)$ of the configuration $\omega_d^\ast$. This problem was first considered by Stolarsky \cite {Sto1975circle,Sto1975} for the Riesz potential functions $f_s(t)$ and certain $s<0$. He proved that $p_{f_s}({\bf x},\omega_2^\ast)$ is minimized on $S^1$ at points of $-\omega_2^\ast$ and maximized at points of $\omega_2^\ast$ for $-2<s<0$ and $-6<s<-4$. However, for $-4<s<-2$, it is minimized at points of $\omega_2^\ast$ and maximized at points of $-\omega_2^\ast$, see \cite {Sto1975circle}. Stolarsky also proved that the absolute maximum of the potential $p_{f_s}({\bf x},\omega_d^\ast)$ on $S^{d-1}$ for $d\geq 3$ and $-2<s<0$ is attained at the points of $\omega_d^\ast$ (cf. \cite {Sto1975}). For potential functions $f$ of form \eqref {Rafailov},
the minimizing property of points from $-\omega_d^\ast$ and the maximizing property of points from $\omega_d^\ast$ for the potential $p_f({\bf x},\omega_d^\ast)$ were established for $s>-2$, $s\neq 0$, and $s<-4$ by Nikolov and Rafailov in \cite {NikRaf2011} ($d=2$) and \cite {NikRaf2013} ($d\geq 3$). At the same time, for $-4<s<-2$, paper \cite {NikRaf2013} showed that at points of $-\omega_d^\ast$ the potential $p_f({\bf x},\omega_d^\ast)$ is maximized while at the points of $\omega_d^\ast$ it is minimized. For $s=0,-2$, and $-4$, the potential $p_{f_s}({\bf x},\omega_d^\ast)$ remains constant over $S^{d-1}$, see \cite [Equation (5.10)]{Sto1975} for the non-trivial case of $s=-4$. 

Papers \cite {Sto1975circle,Sto1975,NikRaf2011,NikRaf2013} also solve the problem of minimizing and maximizing over $S^{d-1}$ the potential (with $f$ as in \eqref {Rafailov}) of other regular point configurations such as the vertices of the regular $N$-gon inscribed in $S^1$ and the vertices of the cube and the cross-polytope inscribed in $S^{d-1}$.
In this paper we establish the following result for the potential of $\omega_d^\ast$.
\begin {theorem}\label {minimum}
Suppose $d\geq 2$ and $f:[0,4]\to (-\infty,\infty]$ is a function continuous on $(0,4]$, differentiable in $(0,4)$ with a concave derivative $f'$ on $(0,4)$ such that $\lim\limits_{t\to 0^+}f(t)=f(0)$. Then the absolute minimum of the potential $p_f({\bf x},\omega_d^\ast)$ over ${\bf x}\in S^{d-1}$ is achieved at every point of the set $-\omega_d^\ast$. Furthermore,
\begin {equation}\label {min}
P_f(\omega_d^\ast,S^{d-1})=\min\limits_{{\bf x}\in S^{d-1}}p_f({\bf x},\omega_d^\ast)=f(4)+d\! \cdot \! f\(2-\frac {2}{d}\).
\end {equation}
If, in addition, $f(0)<\infty$, then the absolute maximum of the potential $p_f({\bf x},\omega_d^\ast)$ over ${\bf x}\in S^{d-1}$ is achieved at every point of $\omega_d^\ast$ with
\begin {equation}\label {max1}
\max\limits_{{\bf x}\in S^{d-1}}p_f({\bf x},\omega_d^\ast)=f(0)+d\! \cdot \! f\(2+\frac {2}{d}\).
\end {equation}
If the concavity of $f'$ is strict on $(0,4)$, then the minimum in \eqref {min} is attained only at the points of $-\omega_d^\ast$, and the maximum in \eqref {max1} is attained only at the points of~$\omega_d^\ast$.
\end {theorem}
In the case $f(0)=\infty$, equality \eqref {max1} holds trivially, since its right-hand side becomes infinite with the points of $\omega_d^\ast$ being the only maximizers in \eqref {max1}. The above theorem is new, for example, for the logarithmic ($d\geq 3$) and Gaussian kernels.

Replacing $f$ with $-f$ in Theorem \ref {minimum} yields the following immediate consequence. 

\begin {corollary}\label {minimum'}
Suppose $d\geq 2$ and $f:[0,4]\to [-\infty,\infty)$ is a function continuous on $(0,4]$, differentiable in $(0,4)$ with a convex derivative $f'$ on $(0,4)$ such that $\lim\limits_{t\to 0^+}f(t)=f(0)$. Then the absolute maximum of the potential $p_f({\bf x},\omega_d^\ast)$ over ${\bf x}\in S^{d-1}$ is achieved at every point of the set $-\omega_d^\ast$. Furthermore,
\begin {equation}\label {max'}
\max\limits_{{\bf x}\in S^{d-1}}p_f({\bf x},\omega_d^\ast)=f(4)+d\! \cdot \! f\(2-\frac {2}{d}\).
\end {equation}
If, in addition, $f(0)>-\infty$, then the absolute minimum of the potential $p_f({\bf x},\omega_d^\ast)$ over ${\bf x}\in S^{d-1}$ is achieved at every point of $\omega_d^\ast$ with
\begin {equation}\label {min'}
P_f(\omega_d^\ast,S^{d-1})=\min\limits_{{\bf x}\in S^{d-1}}p_f({\bf x},\omega_d^\ast)=f(0)+d\! \cdot \! f\(2+\frac {2}{d}\).
\end {equation}

If the convexity of $f'$ is strict on $(0,4)$, then the maximum in \eqref {max'} is attained only at the points of $-\omega_d^\ast$ and the minimum in \eqref {min'} is attained only at the points of~$\omega_d^\ast$.
\end {corollary}

Theorem \ref {minimum} and Corollary \ref {minimum'} explain, in particular, the transition from minimizing property of $-\omega_d^\ast$ to maximizing and vise versa for Riesz $s$-kernel in the above mentioned results from \cite {Sto1975circle} and \cite {NikRaf2013} when $s$ passes through values $-2$ and $-4$. This happens because the derivative $f_s'$ of the potential function changes the direction of its concavity.


%

\section {A simple proof of the optimal covering property of a regular simplex}

Our main results rely on an auxiliary geometric statement, which we prove in this section. It also provides a rather short proof of the optimal covering property of the regular $d$-simplex. The optimal covering problem requires finding positions of $N$ points on a given compact set so that closed balls of the same radius centered at those points cover the set and have their common radius as small as possible.
It is the limiting case as $s\to\infty$ of the maximal polarization problem with respect to the Riesz $s$-kernel. For $N=d+1$ points on $S^{d-1}$ and a point configuration $\omega_d=\{{\bf v}_0,{\bf v}_1,\ldots,{\bf v}_d\}$, we denote
$$
\eta(\omega_d,S^{d-1}):=\max\limits_{{\bf x}\in S^{d-1}}\min\limits_{0\leq i\leq d}\left|{\bf x}-{\bf v}_i\right|.
$$
One is required to find the quantity 
\begin {equation}\label {cov}
\eta_{d+1}(S^{d-1}):=\inf\limits_{\omega_d\subset S^{d-1}}\eta(\omega_d,S^{d-1}),
\end {equation}
and optimal covering point configurations; that is, $(d+1)$-point configurations $\OL\omega_d\subset S^{d-1}$ that attain the infimum on the right-hand side of \eqref {cov}.

The exact solution to the optimal covering problem is known on the unit circle $S^1$, where the set of vertices of a regular $N$-gon inscribed in $S^1$ is optimal for every $N\geq 1$. On the sphere $S^2$, the exact solution is known for $N=1,2,3$ (trivial cases) and for $N=4,6, 12$, see the book by Fejes-T\'oth \cite {Tot1953}. The solution is also known for $N=5, 7, 8, 10$ and $14$ on $S^2$, for $1\leq N\leq 6$ and $N=8$ on $S^3$, and for $1\leq N\leq d+2$ on $S^{d-1}$, $d>4$, see the book by B\"or\"oczky \cite {Bor2004} and references therein.
In particular, the following result is known (see the book by Fejes-T\'oth \cite {Tot1953} for the case $d=3$ and the paper by Galiev \cite {Gal1996} or \cite [Theorem~6.5.1]{Bor2004} for the case $d\geq 4$).
\begin {theorem}\label {covs}
For every $d\geq 2$, a $(d+1)$-point configuration on $S^{d-1}$ provides optimal covering for $S^{d-1}$ if and only if it consists of the vertices of a regular $d$-simplex inscribed in $S^{d-1}$. Furthermore, $\eta_{d+1}(S^{d-1})=\sqrt {2-\frac {2}{d}}$.
\end {theorem}
Known proofs of Theorem \ref {covs} are a bit long compared to the proof of best-packing property of $\omega_d^\ast$. They use spherical geometry (in the case $d=3$), see \cite {Tot1953}, Steiner symmetrization (for any $d\geq 3$), see \cite [Theorem~6.5.1]{Bor2004}, or bounds for a more general problem, see \cite {Gal1996}. We propose a shorter proof that uses the barycentric coordinates of the center ${\bf 0}$ of the sphere $S^{d-1}$ relative to the simplex $T$ formed by the $d+1$ points of the configuration (if ${\bf 0}\notin T$, the proof is elementary).

If a configuration $\omega_d=\{{\bf v}_0,{\bf v}_1\ldots,{\bf v}_d\}\subset S^{d-1}$ is in general position; i.e, not contained in any hyperplane, then the simplex $T$ with the set of vertices $\omega_d$ is non-degenerate and the vertices are affinely independent.
If ${\bf v}$ is any point in $T$, then there is a unique vector $(\beta_0,\beta_1,\ldots,\beta_d)$ of non-negative real numbers such that $\sum_{i=0}^{d}\beta_i=1$ and ${\bf v}=\sum_{i=0}^d\beta_i {\bf v}_i$. The numbers $\beta_i$ are known as the barycentric coordinates of the point ${\bf v}$ relative to the simplex $T$. If ${\bf v}\in {\rm int}\ \! T$, then the barycentric coordinates of ${\bf v}$ are strictly positive, and since their sum is~$1$, each of them is strictly less than $1$.
Let $H_i$, $i=0,1,\ldots,d$, be the hyperplane containing all the points ${\bf v}_j$ except ${\bf v}_i$ and let $L_i$ be the hyperplane parallel to $H_i$ and passing through ${\bf 0}$. Since $\omega_d$ is in general position, we have ${\bf v}_i\notin H_i$. Denote by $H_i^\circ$ the open half-space relative to the hyperplane $H_i$ that contains ${\bf v}_i$ and by $\OL H_i$ the closure of the half-space $H_i^\circ$. Let $r_i$ denote the distance from ${\bf 0}$ to the hyperplane $H_i$ and let $a_i$ denote the distance from the vertex ${\bf v}_i$ to the hyperplane~$L_i$. Denote by $h_i$ the distance from the vertex ${\bf v}_i$ to the hyperplane $H_i$; i.e., the height of $T$. Let also ${\rm conv}(\omega_d)$ be the convex hull of the set $\omega_d$. The proof of Theorem \ref {covs} is essentially contained the following lemma.

\begin {lemma}\label {a_i}
Let ${\bf 0}\in {\rm int}\ \! T$, where $T={\rm conv}(\omega_d)$, and let $b_0,b_1,\ldots,b_d$ be the barycentric coordinates of $\ \!{\bf 0}$ relative to $T$. Then for any index $i$ such that $b_i\leq \frac {1}{d+1}$, we have $r_i\leq 1/d$ and $dr_i\leq a_i$. If $b_i<\frac {1}{d+1}$, then both these inequalities are strict.
\end {lemma}
\begin {proof}
We have $b_i=r_i/h_i\leq \frac {1}{d+1}$.
It is also not difficult to see that $a_i\leq 1$. Then 
\begin {equation}\label {bk}
r_i(d+1)\leq h_i=r_i+a_i\leq r_i+1.
\end {equation}
Consequently, $dr_i\leq 1$ or $r_i\leq 1/d$ and $ dr_i\leq a_i$. If $b_i<\frac {1}{d+1}$, then the first inequality in \eqref {bk} is strict. Consequently, $r_i<1/d$ and $dr_i<a_i$.
\end {proof}

 {\it Proof of Theorem \ref {covs}.}
Let $\omega_d=\{{\bf v}_0,{\bf v}_1,\ldots,{\bf v}_d\}\subset S^{d-1}$ be an arbitrary configuration. If ${\bf 0}\notin {\rm int}\ \! T$, then it is not difficult to see that $\omega_d$ is contained in a hemisphere. Consequently, $\eta(\omega_d,S^{d-1})\geq \sqrt {2}>\sqrt{2-2/d}$.
%
If ${\bf 0}\in {\rm int}\ \! T$, then the key observation is that the average of the barycenric coordinates is $\frac {1}{d+1}$, since they add up to $1$. Then there is an index $k$ such that $b_k\leq \frac {1}{d+1}$. Let $C_k$ be the open spherical cap cut from $S^{d-1}$ by the hyperplane $H_k$ that does not contain the vertex ${\bf v}_k$. Then $C_k\cap\omega_d=\emptyset$. In view of Lemma~\ref {a_i}, the Euclidean radius $R_k$ of $C_k$ satisfies
\begin{equation}\label{eta}
\eta(\omega_d,S^{d-1})\geq R_k=\sqrt{(1-r_k)^2+1-r_k^2}=\sqrt {2-2r_k}\geq \sqrt{2-2/d}=\eta(\omega_d^\ast,S^{d-1}),
\end {equation}
which shows that $\omega_d^\ast$ is best-covering.


Assume that $\eta(\omega_d,S^{d-1})=\sqrt {2-2/d}$. Then ${\bf 0}\in {\rm int}\ \! T$. If for some $k$, $b_k<\frac {1}{d+1}$, then by Lemma~\ref {a_i}, we have $r_k<1/d$ and the second inequality in \eqref {eta} becomes strict contradicting our assumption. Therefore, $b_i\geq \frac {1}{d+1}$ for every $0\leq i\leq d$. Then $b_i=\frac {1}{d+1}$ for all $i$, since $b_i$'s sum to $1$. Then Lemma \ref {a_i} implies that $r_i\leq 1/d$ for all $i$. If $r_k<1/d$ for some $k$, then the second inequality in \eqref {eta} is strict leading again to a contradiction. Thus, $b_i=\frac {1}{d+1}$ and $r_i=1/d$ for all~$i$. Then $\omega_d=\omega_d^\ast$ as the lemma below asserts.
\begin {lemma}\label {1/d+1}
Let ${\bf 0}\in {\rm int}\ \! T$, where $T={\rm conv}(\omega_d)$, and let $b_i=\frac {1}{d+1}$ and $r_i=1/d$ for $0\leq i\leq d$. Then $\omega_d=\omega_d^\ast$.
\end {lemma}
\begin {proof}
Since $a_i\leq 1$, Lemma \ref {a_i} implies that $a_i=1$ for all~$i$. Consequently, ${\bf v}_i\bot L_i$ for all $i$. Hence, ${\bf v}_i\cdot {\bf v}_j=-r_i=-1/d$ for any $i$ and $j\neq i$; that is, $\omega_d=\omega_d^\ast$.
\end {proof}

\section {Proofs of Theorem \ref {minimum} and of Corollary \ref {minimum'}}\label {3}

We start with a basic statement from function theory.
\begin {lemma}\label {AC}
Let a function $g:[a,b]\to \RR$, $a<b$, be continuous on $[a,b]$ and differentiable in $(a,b)$ with a convex derivative $g'$ on $(a,b)$. Then $g$ is absolutely continuous on $[a,b]$.
\end {lemma}
\begin {proof}
Since $g'$ is convex on $(a,b)$, $g'$ is continuous on $(a,b)$. Then the following three cases are possible: (a) $g'$ is non-increasing on $(a,b)$, (b) $g'$ is non-decreasing on $(a,b)$, (c) there is a point $x_0\in (a,b)$ such that $g'$ is non-increasing on $(a,x_0)$ and non-decreasing on $(x_0,b)$. We let $x_0:=b$ in the case (a) and $x_0:=a$ in the case (b). Since $g$ is also continuous on $[a,b]$, using the Mean Value Theorem, we obtain that $g$ is convex on $[x_0,b]$ and concave on $[a,x_0]$ (one of these intervals degenerates into a point in the case (a) or (b)). It is not difficult to see that a convex (concave) monotone continuous function on a finite closed interval is absolutely continuous on that interval. Then each interval $[a,x_0]$ and $[x_0,b]$ can be split into at most two subintervals such that $g$ is absolutely continuous on each. Consequently, $g$ is absolutely continuous on their union, the interval $[a,b]$.
\end {proof}

Next, we establish the inequalities crucial for the proof of Theorem~\ref {minimum}.
\begin {lemma}\label {gt}
Let $d\geq 2$ and $g:[-1,1]\to (-\infty,\infty]$ be a function continuous on $[-1,1)$ and differentiable on $(-1,1)$ with a convex derivative $g'$ on $(-1,1)$ such that $\lim\limits_{t\to 1^-}g(t)=g(1)$. Let $t_0,t_1,\ldots,t_d\in [-1,1]$ be such that
\begin {equation}\label {t_i}
\sum\limits_{i=0}^{d}t_i=0\ \  \ \text {and}\ \ \ \sum\limits_{i=0}^{d}t_i^2=\frac {d+1}{d}.
\end {equation}
Then
\begin {equation}\label {g_q}
g(-1)+d\! \cdot\! g\(\frac {1}{d}\)\leq \sum\limits_{i=0}^{d}g(t_i)\leq g(1)+d\! \cdot \! g\left(-\frac {1}{d}\right).
\end {equation}
If, in addition, $g'$ is strictly convex on $(-1,1)$, then the left inequality in \eqref {g_q} is strict whenever $t_i>-1$, $0\leq i\leq d$, and the right inequality in \eqref {g_q} is strict whenever $t_i<1$, $0\leq i\leq d$.
\end {lemma}
\begin {proof}
We start by establishing the left inequality in \eqref {g_q}.
Assume first that $t_i=-1$ for some $i$. Let $\epsilon_j:=t_j-1/d$, $j\neq i$. Equations \eqref {t_i} imply that $\sum\limits_{j: j\neq i}\epsilon_j=0$ and that
$$
\frac {1}{d}=\sum\limits_{j:j\neq i}t_j^2=\sum\limits_{j:j\neq i}\(\frac {1}{d^2}+\frac {2}{d}\epsilon_j+\epsilon_j^2\)=\frac {1}{d}+\sum\limits_{j:j\neq i}\epsilon_j^2.
$$
Then $\sum\limits_{j:j\neq i}\epsilon_j^2=0$; that is, $\epsilon_j=0$, $j\neq i$. Consequently, $t_j=1/d$, $j\neq i$, and the left inequality in \eqref {g_q} holds as an equality.

Assume next that $t_i>-1$ for all $i$.
If $g(1)=\infty$, then the left inequality in \eqref {g_q} holds trivially whenever one of the $t_i$'s equals $1$. Therefore, in this case we will make an additional assumption that $t_i<1$, for all $i$. 

By Lemma \ref {AC}, the function $g$ is absolutely continuous on any closed subinterval of $[-1,1)$, and, if $g(1)<\infty$, then $g$ is absolutely continuous on $[-1,1]$.
Let $k$ and $\ell\neq k$ be the indices such that $t_k$ is the smallest and $t_{\ell}$ is the second smallest among the numbers $t_i$. Then $t_k>-1$ and by \eqref {t_i}, we have $dt_\ell \leq\sum\limits_{i:i\neq k}t_i<1$. Then $t_\ell<1/d$. Let $A:=\{i\neq k : t_i\geq 1/d\}$ and $B:=\{i\neq k : t_i<1/d\}$ (then $\ell\in B\neq \emptyset$). 
Consider the difference
\begin {equation}\label {p_g}
\begin {split}
\sum\limits_{i=0}^{d}g(t_i)-g(-1)&-d\! \cdot \! g(1/d)=g(t_k)-g(-1)\\
&+\sum\limits_{i\in A}\(g(t_i)-g(1/d)\)-\sum\limits_{i\in B}\(g(1/d)-g(t_i)\)\\
&=\int\limits_{-1}^{t_k}g'(t)\ \!dt+\sum\limits_{i\in A}\int\limits_{1/d}^{t_i}g'(t)\ \! dt-\sum\limits_{i\in B}\ \!\int\limits_{t_i}^{1/d}g'(t)\ \!dt.
\end {split}
\end {equation}
Let $y(t):=\alpha t+\gamma$, where 
$$
\alpha:=\frac {g'(1/d)-g'(t_\ell)}{1/d-t_\ell}\ \ \ \text {and}\ \ \ \gamma=g'\({1}/{d}\)-\frac {\alpha}{d}.
$$
Then $y(t)$ is the polynomial of degree at most one that interpolates $g'(t)$ at $t=t_\ell$ and $t=1/d$. Since $g'$ is convex on $(-1,1)$, we have 
\begin {equation}\label {gconv}
g'(t)\leq y(t),\ \ t\in ( t_\ell,1/d),\ \ \ \text {and}\ \ \ g'(t)\geq y(t),\ \ t\in (-1,t_\ell)\cup (1/d,1).
\end {equation}
Then from \eqref {p_g}, taking into account \eqref {t_i}, we obtain
\begin {equation}\label {pconv}
\begin {split}
\sum\limits_{i=0}^{d}g(t_i)&-g(-1)-d\!\cdot \! g(1/d)\geq \int\limits_{-1}^{t_k}(\alpha t+\gamma)\ \!dt+\sum\limits_{i\in A}\int\limits_{1/d}^{t_i}(\alpha t+\gamma)\ \! dt\\
&-\sum\limits_{i\in B}\ \!\int\limits_{t_i}^{1/d}(\alpha t+\gamma)\ \!dt=\sum\limits_{i=0}^{d}\(\frac {\alpha}{2}t_i^2+\gamma t_i\)-\(\frac {\alpha}{2}-\gamma\)-d\(\frac {\alpha}{2d^2}+\frac {\gamma}{d}\)\\
&=\frac {\alpha}{2}\cdot \frac {d+1}{d}-\frac {\alpha}{2}-\frac {\alpha}{2d}=0.
\end {split}
\end {equation}
Thus, $\sum\limits_{i=0}^{d}g(t_i)\geq g(-1)+d\!\cdot \! g(1/d)$.
Assume that $g'$ is strictly convex on $(-1,1)$ and $t_i>-1$ for all $i$. If $g(1)=\infty$, then the left inequality in \eqref {g_q} holds as a strict inequality whenever $t_i=1$ for some $i$. Therefore, we make an additional assumption that $t_i<1$ for all $i$ when $g(1)=\infty$. 

Both inequalities in \eqref {gconv} are pointwise strict. Then the inequality in \eqref {pconv} is strict and, hence, so is the left inequality in \eqref {g_q}.

Concerning the right inequality in \eqref {g_q}, we observe first that it holds trivially if $g(1)=\infty$. Assume that $g(1)<\infty$ and let $h(t):=-g(-t)$. The function $h$ satisfies the assumptions of the lemma. If numbers $t_0,t_1,\ldots,t_d$ lie in $[-1,1]$ and satisfy \eqref {t_i}, then numbers $\tau_i=-t_i$, $i=0,1,\ldots,d$, also lie in $[-1,1]$ and satisfy \eqref {t_i}. Applying the left inequality in \eqref {g_q} to the function $h$ and the numbers $\tau_i$, we have
\begin {equation}\label {h}
\sum\limits_{i=0}^{d}g(t_i)=-\sum\limits_{i=0}^{d}h(\tau_i)\leq -h(-1)-d\!\cdot \! h(1/d)=g(1)+d\!\cdot \! g(-1/d).
\end {equation}
If $g'$ is strictly convex on $(-1,1)$, then so is $h'$. Since the left inequality in \eqref {g_q} applied to $h$ is strict under the assumption that $\tau_i>-1$ for all $i$, so is the inequality in \eqref {h}. Then the right inequality in \eqref {g_q} is strict whenever $t_i<1$ for all $i$.
\end {proof}

\begin {proof}[Proof of Theorem \ref {minimum}]
Let $g(t):=f(2-2t)$, $t\in [-1,1]$. Then $g$ is continuous on $[-1,1)$ with $\lim\limits_{t\to 1^-}g(t)=g(1)$ and $g'$ is convex on $(-1,1)$. Fix an arbitrary point ${\bf x}\in S^{d-1}$ and let $t_i:={\bf x}\cdot {\bf x}_i$, $0\leq i\leq d$, where ${\bf x}_0,{\bf x}_1,\ldots,{\bf x}_d$ are points in $\omega_d^\ast$. Then $t_i$'s are contained in $[-1,1]$ and, since the center of mass of $\omega_d^\ast$ is at the origin,
$$
\sum\limits_{i=0}^{d}t_i=\sum\limits_{i=0}^{d}{\bf x}\cdot {\bf x}_i={\bf x}\cdot \sum\limits_{i=0}^{d}{\bf x}_i={\bf x}\cdot {\bf 0}=0.
$$
Furthermore, in view of Proposition \ref {squares} from the Appendix and the fact that ${\bf x}\in S^{d-1}$, 
$$
\sum\limits_{i=0}^{d}t_i^2=\sum\limits_{i=0}^{d}({\bf x}\cdot {\bf x}_i)^2=\frac {d+1}{d}.
$$
Thus, relations \eqref {t_i} are satisfied. Then by Lemma \ref {gt},
\begin {equation}\label{minp}
\begin {split}
p_f({\bf x},\omega_d^\ast)&=\sum\limits_{i=0}^{d}f\(2-2{\bf x}\cdot {\bf x}_i\)=\sum\limits_{i=0}^{d}g\({\bf x}\cdot {\bf x}_i\)=\sum\limits_{i=0}^{d}g(t_i)\\
&\geq g(-1)+d\!\cdot \! g(1/d)=f(4)+d\!\cdot \! f(2-2/d)=p_f(-{\bf x}_k,\omega_d^\ast)
\end {split}
\end {equation}
for any $k=0,1,\ldots,d$. This proves relation \eqref {min} and shows that the minimum in \eqref {min} is attained at every point of $-\omega_d^\ast$. If $f'$ is strictly concave on $(0,4)$, then $g'$ is strictly convex on $(-1,1)$. If ${\bf x}\notin -\omega_d^\ast$, then $t_i={\bf x}\cdot {\bf x}_i>-1$ for all $i$, and, by Lemma \ref {gt}, the inequality in \eqref {minp} is strict showing that the minimum in \eqref {min} is not attained at any ${\bf x}\in S^{d-1}\setminus (-\omega_d^\ast)$.

To find points of absolute maximum of the potential $p_f(\cdot , \omega_d^\ast)$ on $S^{d-1}$ ($f(0)<\infty$ implies that $g(1)<\infty$), we use Lemma \ref {gt} and, for every ${\bf x}\in S^{d-1}$, obtain
\begin {equation}\label {uppr}
p_f({\bf x},\omega_d^\ast)=\sum\limits_{i=0}^{d}g(t_i)\leq g(1)+d\!\cdot \! g(-1/d)=f(0)+d\!\cdot \! f(2+2/d)=p_f({\bf x}_k,\omega_d^\ast)
\end {equation}
for any $k=0,1,\ldots,d$. This proves relation \eqref {max1} and shows that the maximum in \eqref {max1} is attained at every point of $\omega_d^\ast$. If $f'$ is strictly concave on $(0,4)$, then $g'$ is strictly convex on $(-1,1)$. If ${\bf x}\notin \omega_d^\ast$, then $t_i={\bf x}\cdot {\bf x}_i<1$ for all $i$, and, by Lemma \ref {gt}, the inequality in \eqref {uppr} is strict showing that the maximum in \eqref {max1} is not attained at any ${\bf x}\in S^{d-1}\setminus \omega_d^\ast$. This completes the proof of Theorem \ref {minimum}. 
\end {proof}

Corollary \ref {minimum'} follows immediately from Theorem \ref {minimum}.



\section {Proof of Theorems \ref {optimalsimplex} and \ref {optimalsimplex'}}

Before proving the main results of this section, we will establish three auxiliary statements.
\begin {lemma}\label {conv}
Let $d\geq 2$ and $f:(0,4]\to \RR$ be a function convex on $(0,4]$. Then the function 
\begin {equation}\label {u(t)}
u(t)=f(2+2dt)+d\!\cdot \! f(2-2t)
\end {equation}
is non-decreasing on $\[0,\frac {1}{d}\]$. If, in addition, $f$ is strictly convex on $(0,4]$, then $u$ is strictly increasing on $\[0,\frac {1}{d}\]$.
\end {lemma}
The case $d=3$ of the above lemma was established in \cite {Su2014}.
\begin {proof}[Proof of Lemma \ref {conv}]
Choose arbitrary $0\leq t_1<t_2\leq \frac {1}{d}$. Then the inequality $u(t_1)\leq u(t_2)$ is equivalent to
\begin {equation}\label{convslope}
\begin {split}
d\!\cdot\! f(2-2t_1)-d\!\cdot \! f(2-2t_2)\leq f(2+2dt_2)-f(2+2dt_1)\\
\frac { f(2-2t_1)- f(2-2t_2)}{2t_2-2t_1}\leq \frac {f(2+2dt_2)-f(2+2dt_1)}{2dt_2-2dt_1},
\end {split}
\end {equation}
which is true in view of convexity of the function $f$, since $$0<2-2t_2<2-2t_1\leq 2+2dt_1<2+2dt_2\leq 4.$$ If $f$ is strictly convex, then a strict inequality holds in \eqref {convslope}, which implies that $u(t_1)<u(t_2)$.
\end {proof}

\begin {lemma}\label {conv'}
Let $d\geq 2$ and $f:[0,4]\to(-\infty,\infty]$ be a function finite and convex on $(0,4]$ such that $f(0)=\lim\limits_{t\to 0^+}f(t)$. Then the function $u(t)$ defined by \eqref {u(t)}
is non-increasing on $\[-\frac {1}{d},0\]$. If, in addition, $f$ is strictly convex on $(0,4]$, then $u$ is strictly decreasing  on $\[-\frac {1}{d},0\]$.
\end {lemma}
\begin {proof}
Choose arbitrary $-\frac {1}{d}< t_1<t_2\leq 0$. Then, by an argument similar to the one in \eqref {convslope}, the inequality $u(t_1)\geq u(t_2)$ is equivalent to 
\begin {equation}\label {decslope}
\frac { f(2-2t_1)- f(2-2t_2)}{2t_2-2t_1}\geq \frac {f(2+2dt_2)-f(2+2dt_1)}{2dt_2-2dt_1},
\end {equation}
which is true in view of convexity of $f$ on $(0,4]$, since $$0<2+2dt_1<2+2dt_2\leq 2-2t_2<2-2t_1<4.$$ Then $u$ is non-increasing on $(-1/d,0]$. If $f$ is strictly convex, then a strict inequality holds in \eqref {decslope}, and, hence, $u$ is strictly decreasing on $(-1/d,0]$. The convexity of $f$ on $(0,4)$ implies its continuity on $(0,4)$. Then $u(-1/d)=
 \lim\limits_{t\to -1/d^+}u(t)$. Consequently, $u$ is monotone on the closed interval $[-1/d,0]$ (and strictly monotone if $f$ is strictly convex).
\end {proof}

In the proof of the three remaining lemmas we will use the notation and definitions introduced between the statements of Theorem \ref {covs} and Lemma \ref {a_i}.

\begin {lemma}\label {semi}
Let $d\geq 2$ and $f:[0,4]\to(-\infty,\infty]$ be a function finite, non-increasing, and convex on $(0,4]$ such that $\lim\limits_{t\to 0^+}f(t)=f(0)$. Suppose $\omega_d\subset S^{d-1}$ is an arbitrary $(d+1)$-point configuration such that ${\bf 0}\notin {\rm int}\ \! T$, where $T={\rm conv}(\omega_d)$. Then 
\begin {equation}\label {z1}
P_f(\omega_d,S^{d-1})\leq \min\left\{f(4)+d\!\cdot\! f\(2-\frac {2}{d}\),f(0)+d\!\cdot\!f\(2+\frac {2}{d}\)\right\}.
\end {equation}
If, in addition, $f$ is strictly convex on $(0,4]$, then the inequality in \eqref {z1} is strict.
\end {lemma}
\begin {proof} 
Since ${\bf 0}\notin {\rm int}\ \! T$, there exists a closed hemi-sphere $H\subset S^{d-1}$ containing $\omega_d:=\{{\bf v}_0,{\bf v}_1,\ldots,{\bf v}_d\}$. Let ${\bf a}\in S^{d-1}$ be such that $H=\{{\bf x}\in S^{d-1} : {\bf x}\cdot {\bf a}\leq 0\}$. By assumption, $f$ is non-increasing on $[0,4]$. Then, since ${\bf v}_i\in H$, $i=0,1,\ldots,d$,
\begin{equation*}
\begin {split}
P_f(\omega_d,S^{d-1})&\leq \sum\limits_{i=0}^{d}f\(\left|{\bf a}-{\bf v}_i\right|^2\)=\sum\limits_{i=0}^{d}f(2-2{\bf a}\cdot {\bf v}_i)\leq (d+1)f(2)=u(0).
\end {split}
\end {equation*}
Taking into account Lemma \ref {conv},  we obtain
\begin {equation}\label {a'}
P_f(\omega_d,S^{d-1})\leq u(1/d)=f(4)+d\!\cdot\! f\(2-\frac {2}{d}\).
\end {equation}
In view of Lemma \ref {conv'}, we also have
\begin {equation}\label {b'}
P_f(\omega_d,S^{d-1})\leq u(-1/d)=f(0)+d\!\cdot\!f\(2+\frac {2}{d}\).
\end {equation}

If $f$ is strictly convex, then by Lemma \ref {conv}, we have $u(0)<u\(1/d\)$ and by Lemma \ref {conv'}, we have $u(0)<u(-1/d)$. Consequently, the  inequality in \eqref {z1} is strict.
\end {proof}

It remains to compare the minimum values over $S^{d-1}$ of potentials of the sets of vertices of an arbitrary simplex whose interior contains ${\bf 0}$ and of a regular simplex. For $i=0,1,\ldots,d$, let ${\bf w}_i$ be the point on $S^{d-1}$ such that the vector from ${\bf 0}$ to ${\bf w}_i$ is perpendicular to the hyperplane $H_i$ and ${\bf w}_i\notin H_i^\circ$.
\begin {lemma}\label {polarizat}
Let $f:[0,4]\to(-\infty,\infty]$ be a non-increasing function finite and convex on $(0,4]$.
Let $d\geq 2$ and $\omega_d\subset S^{d-1}$ be an arbitrary $(d+1)$-point configuration such that ${\bf 0}\in {\rm int}\ \!T$, where $T={\rm conv}\ \! (\omega_d)$. Then
\begin {equation}\label {comp}
P_f(\omega_d,S^{d-1})\leq f(4)+d\!\cdot \! f\(2-\frac {2}{d}\).
\end {equation}
If, in addition, $f$ is strictly convex on $(0,4]$ and $\omega_d$ is not the set of vertices of a regular $d$-simplex inscribed in $S^{d-1}$, then the inequality in \eqref {comp} is strict.
\end {lemma}
We remark that the cases $d=2$ and $d=3$ of Lemma \ref {polarizat} were proved in \cite {Sto1975circle,NikRaf2011,Su2014}.

\begin{proof}[Proof of Lemma \ref {polarizat}]
Recall that $b_0,b_1,\ldots,b_d$ are the barycentric coordinates of ${\bf 0}$ relative to $T$. 
Let $i$ be any index such that $b_i\leq \frac {1}{d+1}$ (since the sum of $b_j$'s is $1$, such an index $i$ exists). 
Observe that ${\bf w}_i\cdot {\bf v}_i=-a_i$ and ${\bf w}_i\cdot {\bf v}_j=r_i$, $j\neq i$. By Lemma \ref {a_i}, we have $-{\bf w}_i\cdot {\bf v}_i\geq dr_i$ and ${\bf w}_i\cdot {\bf v}_j=r_i\leq 1/d$, $j\neq i$. 
Taking into account Lemma \ref {conv}, we now have
\begin {equation}\label {P}
\begin {split}
P_f(\omega_d,S^{d-1})&\leq \sum\limits_{j=0}^{d}f\(\left|{\bf w}_i-{\bf v}_j\right|^2\)=f(2-2{\bf w}_i\cdot {\bf v}_i)+\sum\limits_{j: j\neq i}f\(2-2{\bf w}_i\cdot {\bf v}_j\)\\
&\leq f(2+2dr_i)+d\!\cdot\!f(2-2r_i)\\
&=u(r_i)\leq u\(\frac {1}{d}\)=f(4)+d\!\cdot \! f\(2-\frac {2}{d}\),
\end {split}
\end {equation}
which proves \eqref {comp}.

Assume that $f$ is strictly convex on $(0,4]$ and that $P_f(\omega_d,S^{d-1})=f(4)+d\!\cdot \! f\(2-{2}/{d}\)$. Then $u$ is strictly increasing on $[0,1/d]$ (see Lemma \ref {conv}). If it were that $b_i<\frac {1}{d+1}$, for some $i$, we could write \eqref {P} for that $i$. Lemma \ref {a_i} would imply that $r_i<1/d$ making \eqref {P} strict and contradicting the assumption that equality holds in \eqref {comp}. 
Therefore, $b_i\geq \frac {1}{d+1}$ for every $0\leq i\leq d$, and, consequently, $b_i=\frac {1}{d+1}$ for all $i$, since $b_i$'s sum to $1$. Then for every $i$, we can write \eqref {P} and the assumption that equality holds in \eqref {comp} implies equality throughout \eqref {P}. Furthermore, by Lemma \ref {a_i}, the equality $b_i=\frac {1}{d+1}$ implies that $r_i\leq 1/d$. Since $u$ is strictly increasing on $[0,1/d]$, we must have $r_i=1/d$ for every $i$. Then $\omega_d=\omega^\ast_d$ by Lemma \ref {1/d+1}.
\end {proof}

Finally, we establish the following auxiliary statement.
\begin {lemma}\label {polarizat1}
Let $f:[0,4]\to(-\infty,\infty]$ be a function finite, non-increasing, and convex on $(0,4]$ such that $\lim\limits_{t\to 0^+}f(t)=f(0)$.
Let $d\geq 2$ and $\omega_d\subset S^{d-1}$ be an arbitrary $(d+1)$-point configuration such that ${\bf 0}\in {\rm int}\ \! T$, where $T={\rm conv}\ \!(\omega_d)$. Then
\begin {equation}\label {comp'}
P_f(\omega_d,S^{d-1})\leq f(0)+d\!\cdot \! f\(2+\frac {2}{d}\).
\end {equation}
If, in addition, $f$ is strictly convex on $(0,4]$ and $\omega_d$ is not the set of vertices of a regular $d$-simplex inscribed in $S^{d-1}$, then the inequality in \eqref {comp'} is strict.
\end {lemma}
\begin {proof}[Proof of Lemma \ref {polarizat1}]
If $f(0)=\infty$, then \eqref {comp'} holds trivially as a strict inequality. Therefore, we assume that $f(0)<\infty$ and
let, as before, $b_0,b_1,\ldots,b_d$ be the barycentric coordinates of ${\bf 0}$ relative to $T$.
If for some index $\ell$, we have $r_\ell>1/d$, then
\begin {equation}\label {a0}
\begin {split}
P_f(\omega_d,S^{d-1})&\leq \sum\limits_{j=0}^{d}f\(\left|-{\bf w}_\ell-{\bf v}_j\right|^2\)=f(2+2{\bf w}_\ell\cdot {\bf v}_\ell)+\sum\limits_{j: j\neq \ell}f\(2+2{\bf w}_\ell\cdot {\bf v}_j\)\\
&\leq f(0)+\sum\limits_{j:j\neq \ell}f(2+2r_\ell)\leq f(0)+d\!\cdot \! f\(2+\frac {2}{d}\).
\end {split}
\end {equation}
If $f$ is strictly convex on $(0,4]$, then $f$ decreases strictly on $(0,4]$, and, hence, on~$[0,4]$. Then the last inequality in \eqref {a0} is strict.

Now assume that for every $j=0,1,\ldots,d$, we have $r_j\leq 1/d$.
Let $k$ be any index such that $b_k\geq \frac {1}{d+1}$ (such an index $k$ exists since  $b_i$'s add up to $1$). Since $b_k=r_k/h_k$, we have $r_k(d+1)\geq h_k=r_k+a_k$. Consequently, $a_k\leq dr_k$ and, hence, ${\bf w}_k\cdot {\bf v}_k=-a_k\geq -dr_k$. 
Then taking into account Lemma \ref {conv'} and the fact that $-1/d\leq -r_k<0$, we have
\begin {equation}\label {P'}
\begin {split}
P_f(\omega_d,S^{d-1})&\leq \sum\limits_{j=0}^{d}f\(\left|-{\bf w}_k-{\bf v}_j\right|^2\)=f(2+2{\bf w}_k\cdot {\bf v}_k)+\sum\limits_{j: j\neq k}f\(2+2{\bf w}_k\cdot {\bf v}_j\)\\
&\leq f(2-2dr_k)+d\!\cdot\!f(2+2r_k)\\
&=u(-r_k)\leq u\(-\frac {1}{d}\)=f(0)+d\!\cdot \! f\(2+\frac {2}{d}\),
\end {split}
\end {equation}
which completes the proof of \eqref {comp'}.

Assume that $f$ is strictly convex on $(0,4]$ and that
\begin {equation}\label {**}
P_f(\omega_d,S^{d-1})=f(0)+d\!\cdot \! f\(2+\frac{2}{d}\). 
\end {equation}
Then $f$ is strictly decreasing on $[0,4]$. If $r_\ell>1/d$ for some $\ell$, then the last inequality in \eqref {a0} becomes strict contradicting \eqref {**}. Assume that $r_j\leq 1/d$ for all~$j$. If $b_k>\frac {1}{d+1}$, for some $k$, then in the argument above, we have ${\bf w}_k\cdot {\bf v}_k=-a_k> -dr_k$ and the second inequality in \eqref {P'} becomes strict contradicting \eqref {**}. If $b_i\leq \frac {1}{d+1}$ for all $i$, then $b_i=\frac {1}{d+1}$ for all $i$. If for some $k$, we have $r_k<1/d$, then the last inequality in \eqref {P'} will be strict since $u$ is strictly decreasing (see Lemma \ref {conv'}) contradicting again \eqref {**}. Thus, $b_i=\frac {1}{d+1}$ and $r_i=1/d$ for all $i$. By Lemma \ref {1/d+1}, we now have $\omega_d=\omega_d^\ast$.
 \end {proof}

\begin {proof}[Proof of Theorems \ref {optimalsimplex} and \ref {optimalsimplex'}] Since $f$ is convex on $(0,4)$, it is continuous on $(0,4)$. Since $f$ is non-increasing and convex on $(0,4]$, it is also continuous at $t=4$. Then Theorem \ref {optimalsimplex} follows immediately from Theorem~\ref {minimum} and Lemmas~\ref {semi} and~\ref {polarizat} and Theorem \ref {optimalsimplex'} from Corollary \ref {minimum'} and Lemmas~\ref {semi} and~\ref {polarizat1}.
\end {proof}

\section {Appendix}

In Section \ref {3}, we use the following known result, see \cite [Equation (5.10)]{Sto1975}. For completeness, we give its proof.
\begin {proposition}\label {squares}
Let $\omega_d^\ast=\{{\bf x}_0,{\bf x}_1,\ldots,{\bf x}_d\}$ be the set of vertices of a regular $d$-simplex inscribed in $S^{d-1}$, $d\geq 2$. Then for any point ${\bf x}\in \RR^d$, 
\begin {equation}\label {d+1}
\sum\limits_{i=0}^{d}({\bf x}\cdot {\bf x}_i)^2=\frac {d+1}{d}\left|{\bf x}\right|^2.
\end {equation}
\end {proposition}
\begin {proof}
This proposition can be proved in different ways. We will prove it using induction on dimension. For $d=2$, without loss of generality, we can let ${\bf x}_0=(1,0)$, ${\bf x}_1=\(-1/2,\sqrt {3}/2\)$, and ${\bf x}_2=\(-1/2,-\sqrt {3}/2\)$. Then for every ${\bf x}=(x_1,x_2)\in \RR^2$, we have
$$
\sum\limits_{i=0}^{2}({\bf x}\cdot {\bf x}_i)^2=x_1^2+\frac {1}{4}\(x_1-\sqrt {3}x_2\)^2+\frac {1}{4}\(x_1+\sqrt {3}x_2\)^2=\frac {3}{2}x_1^2+\frac {3}{2}x_2^2=\frac {3}{2}\left|{\bf x}\right|^2.
$$
Assume now that $d>2$ and that equality \eqref {d+1} holds for every ${\bf x}\in \RR^{d-1}$ and prove it for any ${\bf x}\in \RR^d$. Without loss of generality, we can denote ${\bf x}_0=(0,\ldots,0,1)\in \RR^{d}$ and let $H:=\{{\bf x}\in \RR^d : {\bf x}\ \!\bot\ \! {\bf x}_0\}$. Given ${\bf x}\in \RR^d$, let $a\in \RR$ and ${\bf b}\in H$ be such that ${\bf x}=a{\bf x}_0+{\bf b}$. Then $\left|{\bf x}\right|^2=a^2+\left|{\bf b}\right|^2$. Let also ${\bf z}_i\in H$ be such that ${\bf x}_i=-\frac {{\bf x}_0}{d}+{\bf z}_i$, $i=1,\ldots,d$. Observe that ${\bf z}_1,\ldots,{\bf z}_d$ are the vertices of a regular simplex in $H$ inscribed in the sphere of radius $R=\frac{\sqrt {d^2-1}}{d}$ centered at ${\bf 0}$. Then ${\bf x}\cdot {\bf x}_i=\(a{\bf x}_0+{\bf b}\)\cdot \(-\frac {{\bf x}_0}{d}+{\bf z}_i\)=-\frac {a}{d}+{\bf b}\cdot {\bf z}_i$, $i=1,\ldots,d$. Using the induction assumption and the fact that $\sum_{i=1}^{d}{\bf z}_i={\bf 0}$ as the vertices of a regular simplex with the center of mass at ${\bf 0}$, we obtain
\begin {equation*}
\begin {split}
\sum\limits_{i=0}^{d}({\bf x}\cdot {\bf x}_i)^2&=a^2+\sum\limits_{i=1}^{d}\(-\frac {a}{d}+{\bf b}\cdot {\bf z}_i\)^2=\(1+\frac {1}{d}\)a^2-\frac {2a}{d}{\bf b}\cdot \sum\limits_{i=1}^{d}{\bf z}_i\\
&+\sum\limits_{i=1}^{d}({\bf b}\cdot {\bf z}_i)^2=\(1+\frac {1}{d}\)a^2+R^2\sum\limits_{i=1}^{d}\({\bf b}\cdot \frac {{\bf z}_i}{R}\)^2\\
&=\frac {d+1}{d}a^2+\frac {dR^2}{d-1}\left|{\bf b}\right|^2=\frac {d+1}{d}(a^2+\left|{\bf b}\right|^2)=\frac {d+1}{d}\left|{\bf x}\right|^2,
\end {split}
\end {equation*}
and \eqref {d+1} follows by induction.
\end {proof}

\end {document}